\documentclass[a4paper,twoside,12pt]{amsart}
\RequirePackage{amsmath,amssymb,xypic}
\theoremstyle{plain}
\newtheorem{theorem}{Theorem}[section]
\newtheorem{corollary}[theorem]{Corollary}
\newtheorem{proposition}[theorem]{Proposition}
\newtheorem{lemma}[theorem]{Lemma}

\theoremstyle{definition}
\newtheorem{definition}[theorem]{Definition}
\newtheorem{example}[theorem]{Example}
\newtheorem{remark}[theorem]{Remark}

\numberwithin{equation}{section}

\newcommand{\Union}{\bigcup\limits}

\newcommand{\C}{\mathbb{C}}
\newcommand{\N}{\mathbb{N}}

\newcommand{\Z}{\mathbb{Z}}

\newcommand{\op}{\mathrm{op}}

\DeclareMathOperator{\id}{id}

\newcommand{\DSum}{\bigoplus}

\DeclareMathOperator{\coker}{coker}

\renewcommand{\to}[1][]{\xrightarrow[]{#1}}

\newcommand{\isoto}[1][]{\xrightarrow[#1]{\sim}}

\newcommand{\Endo}[1][]{\mathrm{End}_{\raise1.5ex\hbox to.1em{}#1}}

\newcommand{\Hom}[1][]{\mathrm{Hom}_{\raise1.5ex\hbox to.1em{}#1}}

\newcommand{\RHom}[1][]{\mathrm{RHom}_{\raise1.5ex\hbox to.1em{}#1}}

\newcommand{\Ext}[2][]{\mathrm{Ext}_{\raise1.5ex\hbox to.1em{}#1}^{#2}}

\newcommand{\THom}[1][]{\mathrm{THom}_{\raise1.5ex\hbox to.1em{}#1}}

\newcommand{\Tens}[1][]{\mathbin{\otimes_{\raise1.5ex\hbox to-.1em{}#1}}}

\newcommand{\LTens}[1][]{\mathbin{\otimes_{\raise1.5ex\hbox to-.1em{}#1}^{L}}}

\newcommand{\Tor}[2][]{\mathrm{Tor}^{\raise1.5ex\hbox to.1em{}#1}_{#2}}

\def\sha{\mathcal{A}}
\def\shb{\mathcal{B}}

\def\she{\mathcal{E}}

\def\shg{\mathcal{G}}

\def\shl{\mathcal{L}}

\def\sho{\mathcal{O}}
\def\shp{\mathcal{P}}

\def\shr{\mathcal{R}}
\def\shs{\mathcal{S}}

\newcommand{\shendo}[1][]{{\mathcal{E}nd}_{\raise1.5ex\hbox to.1em{}#1}}

\newcommand{\shaut}[1][]{{\mathcal{A}ut}_{\raise1.5ex\hbox to.1em{}#1}}

\renewcommand{\hom}[1][]{{\mathcal{H}om}_{\raise1.5ex\hbox to.1em{}#1}}

\newcommand{\rhom}[1][]{{R\mathcal{H}om}_{\raise1.5ex\hbox to.1em{}#1}}

\newcommand{\ext}[2][]{{\mathcal{E}xt}_{\raise1.5ex\hbox to.1em{}#1}^{#2}}

\newcommand{\thom}[1][]{{T\mathcal{H}om}_{\raise1.5ex\hbox to.1em{}#1}}

\newcommand{\tens}[1][]{\mathbin{\otimes_{\raise1.5ex\hbox to-.1em{}#1}}}

\newcommand{\ltens}[1][]{\mathbin{\otimes_{\raise1.5ex\hbox to-.1em{}#1}^{L}}}

\newcommand{\tor}[2][]{{\mathcal{T}or}^{\raise1.5ex\hbox to.1em{}#1}_{#2}}

\newcommand{\opb}[1]{#1^{-1}}

\newcommand{\GHom}[1][]{\mathrm{GHom}_{\raise1.5ex\hbox to.1em{}#1}}

\newcommand{\GExt}[2][]{\mathrm{GExt}_{\raise1.5ex\hbox to.1em{}#1}^{#2}}

\newcommand{\FHom}[1][]{\mathrm{FHom}_{\raise1.5ex\hbox to.1em{}#1}}

\newcommand{\ghom}[1][]{{\mathcal{GH}om}_{\raise1.5ex\hbox to.1em{}#1}}

\newcommand{\gext}[2][]{{\mathcal{GE}xt}_{\raise1.5ex\hbox to.1em{}#1}^{#2}}

\newcommand{\fhom}[1][]{{\mathcal{FH}om}_{\raise1.5ex\hbox to.1em{}#1}}

\newcommand{\gr}{\mathop{\mathcal{G}r}\nolimits}
\newcommand{\Gr}{\mathop{\mathrm{Gr}}\nolimits}

\newcommand{\tenstop}[1][]{\mathbin{\hat{\otimes}_{\raise1.5ex\hbox to-.1em{}#1}}}

\newcommand{\homtop}[1][]{\mathcal{L}_{\raise1.5ex\hbox to.1em{}#1}}

\newcommand{\Homtop}[1][]{\mathrm{L}_{\raise1.5ex\hbox to.1em{}#1}}

\newcommand{\E}{\mathcal{E}}

\renewcommand{\O}{\mathcal{O}}

\def\absdoim#1{\underline{#1}_*}
\def\reldoim[#1]#2{\underline{#2}_{|{#1}*}}
\def\doim{\@ifnextchar [{\reldoim}{\absdoim}}

\def\absdeim#1{\underline{#1}_*}
\def\reldeim[#1]#2{\underline{#2}_{|{#1}*}}
\def\deim{\@ifnextchar [{\reldeim}{\absdeim}}

\def\absdopb#1{\underline{#1}^{-1}}
\def\reldopb[#1]#2{\underline{#2}_{|{#1}}^{-1}}
\def\dopb{\@ifnextchar [{\reldopb}{\absdopb}}

\def\absboim#1{\underline{\underline{#1}}_*}
\def\relboim[#1]#2{\underline{\underline{#2}}_{|{#1}*}}
\def\boim{\@ifnextchar [{\relboim}{\absboim}}

\def\absbeim#1{\underline{\underline{#1}}_*}
\def\relbeim[#1]#2{\underline{\underline{#2}}_{|{#1}*}}
\def\beim{\@ifnextchar [{\relbeim}{\absbeim}}

\def\absbopb#1{\underline{\underline{#1}}^*}
\def\relbopb[#1]#2{\underline{\underline{#2}}_{|{#1}}^*}
\def\bopb{\@ifnextchar [{\relbopb}{\absbopb}}

\newcommand{\ad}{\operatorname{ad}}

\newcommand{\catc}{\mathsf{C}}
\newcommand{\catGr}[1]{\mathsf{Gr}(#1)}
\newcommand{\catMod}{\mathsf{Mod}}
\newcommand{\catFun}[1][]{\mathsf{Hom}_{#1}}

\newcommand{\stkMod}[1][]{\mathfrak{Mod}_{#1}}
\newcommand{\stkAut}[1][]{\mathfrak{Aut}_{#1}}
\newcommand{\stkFun}[1][]{\mathfrak{Hom}_{#1}}

\newcommand{\stkGr}{\mathfrak{Gr}}

\newcommand{\stka}{\mathfrak{A}}
\newcommand{\stke}{\mathfrak{E}}
\newcommand{\stkg}{\mathfrak{G}}

\newcommand{\stks}{\mathfrak{S}}

\newcommand{\astk}[1]{#1^+}

\newcommand{\approxto}[1][]{\xrightarrow[#1]{\approx}}

\renewcommand{\E}{\she^{\sqrt v}}

\begin{document}

\title[Uniqueness of quantization of complex contact manifolds]
{Uniqueness of quantization of complex contact manifolds}

\author[P. Polesello]{Pietro Polesello}
\address{Universit{\`a} di Padova\\
Dipartimento di Matematica Pura ed Applicata\\
via G. Belzoni, 7\\ 
35131 Padova, Italy}
\email{pietro@math.unipd.it}

\thanks{The author had the occasion of visiting Paris VI University 
during the preparation of this paper. Their hospitality is gratefully 
acknowledged.}

\date{}


\begin{abstract}
Using the language of algebroid stacks, we will show that Kashiwara's quantization 
of a complex contact manifold is unique.\\
\end{abstract}

\maketitle

\section*{Introduction} 
Let $M$ be a complex manifold and $P^*M$ its projective cotangent bundle, 
endowed with the canonical contact structure. Let $\she_M$ be the sheaf of 
algebras on $P^*M$ of microdifferential operators. Recall that the order of the 
operators defines a filtration on $\she_M$ such that its associated graded 
algebra is isomorphic to $\DSum_{m\in\Z} \O_{P^*M}(m)$. (Here $ \O_{P^*M}(m)$ 
denotes the $m$-th tensor power of the dual of the tautological bundle 
$\O_{P^*M}(-1)$). The product on $\she_M$ is given by the Leibniz rule and it 
is compatible with the Jacobi structure on $\O_{P^*M}(1)$ induced by the 
Poisson bracket on $T^*X$. Hence this algebra provides a quantization of $P^*M$.
Any filtered sheaf of algebras which has $\DSum_{m\in\Z} \O_{P^*M}(m)$ as 
graded algebra and which is locally isomorphic to $\she_M$ gives another 
quantization of $P^*M$. We call such an object an $\she$-algebra.

On a complex contact manifold $Y$ there may not exist an $\she$-algebra, that 
is, a filtered sheaf of algebras which has $\DSum_{m\in\Z} \shl^{\otimes m}$ 
as graded algebra (here $\shl$ is the line bundle associated to the contact structure) 
and which is locally isomorphic to $\opb i \she_M$, for any contact 
local chart $i\colon Y\supset U \to P^*M$. However, 
Kashiwara~\cite{Kashiwara1996} proved that the stack (sheaf of categories) 
$\stkMod(\she;Y)$ of modules over these locally defined sheaves of algebras is always defined. 

This quantization of $Y$ is better understood using the language of 
algebroid stacks~\cite{Kontsevich2001,D'Agnolo-Polesello2005}. In the spirit of 
Mitchell's notion of ``algebra with several objects''~\cite{Mitchell}, we may say 
that an algebroid stack is a sheaf of algebras with locally several (locally) 
isomorphic objects. With this notion at hand, we may reformulate 
Kashiwara's result by saying that there exists an algebroid stack $\stke_Y$ over $Y$ 
such that $\stkMod(\she;Y)$ is equivalent to the stack of $\stke_Y$-modules. 
The advantage is that $\stke_Y$ has now similar properties to that of an 
$\she$-algebra: it is filtered with $\DSum_{m\in\Z} \shl^{\otimes m}$ 
as associated graded trivial algebroid and locally 
equivalent to the trivial algebroid $\opb i \she_M$, for any contact local 
chart $i\colon Y\supset U \to P^*M$.

The purpose of this paper is to show that the algebroid stack $\stke_Y$ is the unique 
quantization of $Y$ endowed with an anti-involution. 
(We note here that the complex symplectic case behaves differently. Reference is made to
~\cite{Polesello}.)

The paper is organized as follows: in Section~\ref{section:E-alg} we recall
the definition of microdifferential operator and that of $\she$-algebra on 
$P^*M$. In Section~\ref{section:filt-grad} we give the main definitions and 
properties of filtered and graded stacks. In Section~\ref{section:uniq} we 
prove the uniqueness of the algebroid stack $\stke_Y$ (Theorem~\ref{th:unicity}).
In Appendix A we recall some basic facts about stacks of 2-groups, which are 
needed in the proof of Theorem~\ref{th:unicity}.

\medskip
\noindent
{\bf Acknowledgements}
We wish to thank Masaki Kashiwara for valuable comments and insights. We also 
would like to thank Louis Boutet de Monvel, Matthieu Carette and Pierre Schapira 
for useful discussions.

\medskip
\noindent
{\bf Notations and conventions}
All the filtrations are intended to be over $\Z$, increasing and exhaustive. 
If $A$ (resp. $\sha$) is a filtered algebra (resp. sheaf of filtered algebras),
we will denote by $\Gr (A)$ (resp. $\gr(\sha)$) its associated graded algebra
(resp. sheaf of graded algebras), and by $\Gr_0 (A)$ (resp. $\gr_0(\sha)$) the 
algebra (resp. sheaf of algebras) of homogeneous elements of degree $0$. We will use 
similar notations for morphisms.

If $\sha$ is a sheaf of algebras, we will denote by $\sha^\times$ the sheaf of 
groups of its invertible elements and, for each section $a\in \sha^\times$, 
by $\ad(a)\colon \sha \to \sha$ the algebra isomorphism $b\mapsto aba^{-1}$.

We will use the upper index $\op$ to denote opposite structure, when referring 
either to (sheaves of) algebras, to categories or to stacks.

\section{$\she$-algebras on $P^*M$}\label{section:E-alg}

We recall here the basic properties of the algebra of microdifferential operators.
References are made to~\cite{S-K-K,Kashiwara1986,Kashiwara2000} (see 
also~\cite{Schapira} for an exposition).

\medskip

Let $M$ be a complex manifold, and $\pi\colon P^*M\to M$ its projective 
cotangent bundle. Denote by $\she_M$ the sheaf of microdifferential operators, 
considered as a sheaf on $P^*M$. In a local coordinate system $(x)$ on $M$, 
with associated local coordinates $(x;[\xi])$ on $P^*M$, a 
microdifferential operator $P$ of order $m$ defined on an open subset $U$ of 
$P^*M$ has a total symbol
$$
\sigma_{\mathrm{tot}}(P)=\sum_{j=-\infty}^m p_j(x;\xi),
$$
where the $p_j$'s are sections on $U$ of $\O_{P^*M}(j)$, subject to the 
estimates
\begin{equation}\label{eq:estmicrod}
     \left\{ \begin{array}{l}
     \mbox{for any compact subset $K$ of $U$ there exists a constant}\\
     \mbox{$C_K>0$ such that for all $j<0$,}
     \sup\limits_{K}\vert p_{j}\vert \leq C_K^{-j}(-j)!.
     \end{array}\right.
\end{equation}
(Here $ \O_{P^*M}(j)$ denotes the $j$-th tensor power of the dual 
of the tautological bundle $\O_{P^*M}(-1)$).
The product structure on $\she_M$ is given by the Leibniz formula: if $Q$ is 
another microdifferential operator defined on $U$ of total symbol 
$\sigma_{\mathrm{tot}}(Q)$, then 
$$\sigma_{\mathrm{tot}}(P\circ Q)=\sum_{\alpha\in\N^n} \frac{1}
{\alpha !} \partial^{\alpha}_{\xi}\sigma_{\mathrm{tot}}(P)
\partial^{\alpha}_x\sigma_{\mathrm{tot}}(Q).
$$

Recall that the center of $\she_M$ is the constant sheaf $\C_{P^*M}$ and that
$\she_M$ is filtered. One denotes by $\she_M(m)$ the sheaf of 
operators of order less than or equal to $m$ and by 
$$
\sigma_m(\cdot)\colon \she_M(m)\to \she_M(m)/\she_M(m-1)
\simeq \sho_{P^*M}(m)
$$ 
the symbol map of order $m$, which does not depend on the local 
coordinate system on $P^*M$. If $\sigma_m(P)$ is not
identically zero, then one says that $P$ has order $m$ and
$\sigma_m(P)$ is called the principal symbol of $P$. 
In particular, an element $P$ in $\she_M$ is invertible if and only if
its principal symbol is nowhere vanishing.

\begin{remark}
The algebra $\she_M$ is a quantization of $P^*M$ in the following 
sense. Denote by $\shs_{P^*M}$ the sheaf on $P^*M$ whose local 
sections are symbols, that is, series $\sum_{j=-\infty}^m f_j$ with  
$f_j\in\O_{P^*M}(j)$ satisfying the estimates~\eqref{eq:estmicrod}.
It is endowed with a natural filtration and the associated graded 
sheaf is $\DSum_{m\in\Z} \O_{P^*M}(m)$.
Then the sheaf $\she_M$ is locally isomorphic to $\shs_{P^*M}$ as filtered 
$\C_{P^*M}$-modules (via the total symbol), and the product on 
$\shs_{P^*M}$ induced by the Leibniz rule is unitary, associative and 
compatible with the Jacobi structure on $\O_{P^*M}(1)$, that is, the 
following diagram commutes
\begin{equation}
\label{eq:Jacobi}
\xymatrix@C3.5em{ \she_M(1)\times \she_M(1) \ar[r]^-{[\cdot,\cdot]} 
\ar[d]_-{\sigma_1\times\sigma_1} & 
\she_M(1)\ar[d]^-{\sigma_1} \\
\sho_{P^*M}(1) \times \sho_{P^*M}(1) \ar[r]_-{\{\cdot,\cdot\}} &
\sho_{P^*M}(1) .}
\end{equation}
(Here $[\cdot,\cdot]$ denotes the commutator and $\{\cdot,\cdot\}$ is induced 
by the Poisson braket on $T^*M$.)
\end{remark}


Let $\Omega_M$ be the canonical line bundle on $M$, that is, the sheaf of forms of 
top degree. Recall that each locally defined volume form $\theta\in\Omega_M$
gives rise to a local isomorphism $*_{\theta}\colon\she_M^{\op}\isoto \she_M$, 
which sends an operator $P$ to its formal adjoint $P^{*_{\theta}}$ with 
respect to $\theta$.
In a local coordinate system $(x)$ satisfying $\theta = dx$, 
with associated local coordinates $(x;[\xi])$, one has
$$
\sigma_{\mathrm{tot}}(P^{*_{\theta}})=\sum_{\alpha\in\N^n} \frac{(-1)^{|\alpha |}}
{\alpha !} \partial^{\alpha}_{\xi}\partial^{\alpha}_x\sigma_{\mathrm{tot}}(P)(x,-\xi).
$$
Twisting $\she_M$ by $\Omega_M$, one then gets a globally 
defined isomorphism of algebras
$$
\she_M^{\op} \isoto \opb {\pi}\Omega_M\tens\she_M \tens \opb{\pi}
\Omega_M^{\tens -1}
\qquad P\mapsto \theta\tens P^{*_{\theta}}\tens \theta^{\tens -1}, 
$$
which does not depend on the choice of the volume form. 
(Here the tensor product is over $\opb{\pi}\O_M$ and $\theta^{\tens -1}$ 
denotes the unique section of $\Omega_M^{\tens -1}$, the dual of $\Omega_M$, 
satisfying $\theta^{\tens -1}\tens \theta =1$.) 
This leads to replace the algebra $\she_M$ by its twisted version 
by half-forms\footnote{Recall 
that the sections of $\E_M$ are locally defined by $\theta^{\tens 1/2}\tens 
P\tens \theta^{\tens -1/2}$ for a volume form $\theta$ and a microdifferential 
operator $P$, with the equivalence relation $\theta_1^{\tens 1/2}\tens P_1\tens 
\theta_1^{\tens -1/2} = \theta_2^{\tens 1/2}\tens P_2\tens 
\theta_2^{\tens -1/2}$ if and only if $P_2 = (\theta_1/\theta_2)^{1/2} P_1 
(\theta_1/\theta_2)^{-1/2}$.}
$$
\E_M =
\opb {\pi} \Omega_M^{\tens 1/2}\tens\she_M\tens\opb{\pi}\Omega_M^{\tens -1/2}.
$$

The sheaf $\E_M$ is a $\C$-algebra on $P^*M$ locally isomorphic to $\she_M$, 
and it has the following properties:
\begin{itemize} 
\item[(i)] it is filtered;
\item[(ii)] there is an isomorphism of graded algebras 
$$\sigma\colon\gr(\E_M)\isoto \DSum_{m\in\Z} \O_{P^*M}(m);$$
\item[(iii)] it is endowed with an anti-involution $*$, {\em i.e.} 
an isomorphism of algebras $$*\colon (\E_M)^{\op}\isoto\E_M \quad 
\mbox{such that $*^2=\id$.}$$
\end{itemize}
Moreover, these data are compatible: the anti-involution respects the filtration 
and $\gr_0(*)$ is the identity. Note that one has 
$\sigma_m(P^*) = (-1)^m\sigma_m(P)$ for all $P\in \E_M(m)$.


\medskip

This suggests the following
\begin{definition}\label{definition:E-alg}
An $\she$-algebra with anti-involution on $P^*M$, an ($\she,*$)-algebra for short,
is a sheaf of $\C$-algebras $\sha$ 
together with
\begin{itemize}
\item[(i)] a filtration $\{F_m\sha\}_{m\in\Z}$; 
\item[(ii)] an isomorphism of graded algebras $\nu\colon\gr(\sha)\isoto
 \DSum_{m\in\Z} \O_{P^*M}(m)$; 
\item[(iii)] an anti-involution $\iota$;
\end{itemize} 
such that the triplet $(\sha,\nu,\iota)$ is locally isomorphic to 
$(\E_M, \sigma, *)$. 

A morphism of ($\she,*$)-algebras is a $\C$-algebra morphism 
compatible with the structures (i), (ii) and (iii).
\end{definition}

\begin{remark}
Definition~\ref{definition:E-alg} is adapted from~\cite{BoutetdeMonvel1999}. 
See also~\cite{BoutetdeMonvel2002} for similar definitions in the framework 
of real manifolds, where $\she$-algebras are replaced by Toeplitz algebras. 
\end{remark}

By definition, a morphism $\varphi\colon \sha_1\to \sha_2$ of 
($\she,*$)-algebras is a $\C$-algebra morphism commuting with 
the anti-involutions, mapping $F_m\sha_1$ to $F_m\sha_2$ in such a way that 
$\nu^2_m(\varphi(P)) = \nu^1_m(P)$ for all $P\in F_m\sha_1$. 
(Here  $\nu^i_m$ denotes the symbol map 
$F_m(\sha_i)\to F_m(\sha_i)/F_{m-1}(\sha_i)\simeq \sho_{P^*M}(m)$ of order $m$,
for $i=1,2$.) It follows that any ($\she,*$)-algebra provides a 
quantization of $P^*M$. 

\begin{example}
Let $f\colon P^*M \to P^*M$ be a contact transformation. Then $\opb f \E_M$ 
inherits an anti-involution and a filtration from $\E_M$, in such a way that the 
associated graded algebra is isomorphic (via $f$) to $\DSum_{m\in\Z} \O_{P^*M}(m)$. 
By a result of~\cite{S-K-K}, for each $p\in P^*M$ there exist an open neighborhood $V$
of $p$ and an isomorphism\footnote{This isomorphism is called Quantized 
Contact Transformation over $f$.} $\opb f \E_M|_V\isoto \E_M|_V$ respecting 
(i), (ii) and (iii). 
It follows that $\opb f \E_M$ is an ($\she,*$)-algebra.
\end{example}

Denote by $\shaut[\she,*](\E_M)$ the group of automorphisms of $\E_M$ as an
($\she,*$)-algebra and set
$$\she^{\sqrt v, *}_M = \{P\in\E; \mbox{ $P$ has order 0, $\sigma_0(P)=1$ 
and  $PP^*=1$}\}\subset (\she^{\sqrt v}_M)^\times.$$

\begin{lemma}[cf \cite{Kashiwara1996}]\label{lemma:key}
The assignment $P\mapsto \ad(P)$ defines an isomorphism of 
sheaves of groups on $P^*M$
\begin{equation*}
\ad\colon\she^{\sqrt v, *}_M \to  \shaut[\she,*](\E_M).
\end{equation*}
\end{lemma}

The set of isomorphism classes of ($\she,*$)-algebras 
on $P^*M$ is in bijection with $H^1(P^*M; \shaut[\she,*](\E_M))$.
\begin{corollary}
The ($\she,*$)-algebras on $P^*M$ are classified by the pointed set
$H^1(P^*M;\she^{\sqrt v, *}_M)$.
\end{corollary}

Adapting a result in~\cite{BoutetdeMonvel1999}, we get the following

\begin{proposition}
If $\dim M \geq 3$, then $\E_M$ is the unique, up to isomorphism, 
($\she,*$)-algebra on $P^*M$.  
\end{proposition}

\section{Filtered and graded stacks}\label{section:filt-grad}

According to Kashiwara's result~\cite{Kashiwara1996}, one has to replace 
sheaves by stacks in order to quantize a contact complex manifold. 
It becomes thus necessary to define the notions of filtration and graduation 
of a stack. 
We start here by recalling what a filtered (resp. graded) category is and how 
to associate a graded category to a filtered one. Then we stackify these 
definitions.
We assume that the reader is familiar with the basic notions from the theory 
of stacks which are, roughly speaking, sheaves of categories. (The classical 
reference is~\cite{Giraud1971}, and a short presentation is given {\em e.g.}
in~\cite{Kashiwara1996,D'Agnolo-Polesello2003}.)

\medskip
Let $R$ be a commutative ring. 
\begin{definition}\label{def:filtered}
     An $R$-category\footnote{An $R$-category is a category whose 
     sets of morphisms are endowed with an $R$-module structure, so that 
     the composition is bilinear. An $R$-functor is a functor between 
     $R$-categories which is linear at the level of morphisms.} 
     $\catc$ is filtered (resp. graded) if the following properties are 
     satisfied:
     \begin{itemize}
     \item[$\bullet$] for any objects $P,Q\in\catc$,  the $R$-module 
     $\Hom[\catc](P,Q)$ is filtered (resp. graded);
     \item[$\bullet$] for any $P,Q,R\in \catc$ and any morphisms $f$ in 
     $F_m\Hom[\catc](Q,R)$ (resp. in $G_m\Hom[\catc](Q,R)$) and $g$ in 
     $F_n\Hom[\catc](P,Q)$ (resp. in $G_n\Hom[\catc](P,Q)$), the composition 
     $f\circ g$ is in $F_{m+n}\Hom[\catc](P,R)$ (resp. in 
     $G_{m+n}\Hom[\catc](P,R)$);
     \item[$\bullet$] for each $P\in \catc$, the identity morphism $\id_{P}$ 
     is in $F_0\Hom[\catc](P,P)$ (resp. in $G_0\Hom[\catc](P,P)$). 
     \end{itemize}
     
     \noindent
     An $R$-functor $\Phi\colon \catc \to \catc'$ between filtered 
     (resp. graded) categories is filtered (resp. graded) if for any objects 
     $P,Q\in\catc$, the $R$-module morphism $\Hom[\catc](P,Q)\to 
     \Hom[\catc'](\Phi(P),\Phi(Q))$ is filtered (resp. graded).
     
     \noindent
     A natural transformation $\alpha=(\alpha_P)_{P\in \catc}
     \colon\Phi\Rightarrow \Phi'$ between filtered (resp. graded) functors is 
     filtered (resp. graded) if  for each $P\in\catc$ the morphism $\alpha_P$ is in 
     $F_0 \Hom[\catc'](\Phi(P),\Phi'(P))$  (resp. in  
     $G_0 \Hom[\catc'](\Phi(P),\Phi'(P))$).
\end{definition}

If $\Phi, \Phi'\colon \catc \to \catc'$ are filtered $R$-functors, one 
defines $F_n\Hom (\Phi, \Phi')$ as the set of those 
natural transformations of functors 
$\alpha = (\alpha_P)_{P\in\catc}\colon \Phi \Rightarrow \Phi'$ with 
$\alpha_P\in F_n\Hom[\catc'] (\Phi(P), \Phi'(P))$ for each $P\in \catc$.
The category $\catFun[F](\catc,\catc')$ thereby obtained is then filtered. 
A similar remark holds also for the category of graded functors between 
graded categories. 

To any filtered $R$-category $\catc$ there is an associated graded 
$R$-category $\catGr \catc$, whose objects are the same of those of $\catc$ 
and for any objects $P,Q$ the set of morphisms is defined by 
$\Hom[\catGr \catc](P,Q) = \Gr(\Hom[\catc](P,Q))$. Similarly, to any 
filtered functor $\Phi$, one associates a graded functor $\Gr(\Phi)$.
In this way, we get a 
functor $\catGr \cdot$ from filtered $R$-categories to graded ones.
Note that, if $\Phi, \Phi'$ are filtered functors, then 
there is a natural injective morphism 
$\Gr_n(\Hom (\Phi, \Phi'))\to G_n\Hom (\catGr \Phi, \catGr {\Phi'})$ for each $n\in\Z$.
We will denote by $\catGr\alpha$ the graded natural transformation 
$\catGr {\Phi} \Rightarrow \catGr {\Phi'}$ associated to the filtered natural 
transformation $\alpha\colon\Phi\Rightarrow \Phi'$ via the previous 
morphism for $n=0$.\footnote{It follows that $\catGr \cdot$ defines a 
2-functor from the 2-category of filtered categories, filtered functors and 
filtered natural transformation to that of graded categories, 
graded functors and graded natural transformations.}

When we restrict to morphisms homogeneous of degree $0$, we will 
use the notation $\mathsf{Gr}_0$ (for categories, functors and natural transformations).

\medskip

Following the presentation in~\cite{D'Agnolo-Polesello2005}, recall that there 
is a fully faithful functor $(\cdot)^+$ from filtered (resp. graded) 
$R$-algebras to filtered (resp. graded) $R$-categories, which sends an
$R$-algebra $A$ to the category $\astk A$ with a single 
object $\bullet$  and $\Endo(\bullet)=A$ as set of morphisms. 
If $f,g\colon A\to B$ are filtered (resp. graded) $R$-algebra morphisms, then 
the filtered (resp. graded) natural transformations 
$\astk f \Rightarrow \astk g$ correspond to elements $b$ in $F_0 B$ (resp. 
in $G_0 B$) such that $b f(a) = g(a) b$ for any $a\in A$.
Clearly, for any filtered $R$-algebra $A$ one has 
$\catGr {\astk A}=\astk{\Gr (A)}$. Moreover, if $f,g\colon A\to B$ are filtered
$R$-algebra morphisms and $b$ in $F_0 B$ defines a filtered natural 
transformation $\alpha\colon\astk f \Rightarrow \astk g$, then the 
graded natural transformation $\catGr \alpha\colon\astk {\Gr (f)} 
\Rightarrow \astk {\Gr (g)}$ is given by the image of $b$ 
via the natural map $F_0 B \to\Gr_0 B$.

Let $A$ be a filtered $R$-algebra. For any filtered (left) $A$-modules $M$ 
and $N$, let $F_m\Hom[FA](M,N)$ be the set of those $A$-module morphisms 
$M \to N$ which send $F_n M$ to $F_{n+m}N$. We denote by
$\catMod_F(A)$ the filtered $R$-category so obtained. 
One easily checks that it is equivalent to the category 
$\catFun[F](\astk A,\catMod_F(R))$, where $R$ is equipped with the trivial 
filtration, and that the Yoneda embedding
$$
\astk A \to \catFun[F]((\astk A)^\op, \catMod_F(R)) \approx 
\catMod_F(A^\op)
$$
identifies $\astk A$ with the full subcategory of filtered right $A$-modules 
which are isomorphic to $A$. Everything remains true replacing filtered 
algebras and categories by graded ones.

\medskip

Let $X$ be a topological space, and $\shr$ a sheaf of commutative rings.

As for categories, there are natural notions of filtered (resp. graded) 
$\shr$-stack, of filtered (resp. graded) $\shr$-functor between 
filtered (resp. graded) stacks and of filtered (resp. graded) natural 
transformations between filtered (resp. graded) functors.

 
As above, we denote by $\astk{(\cdot)}$ the faithful and locally full 
functor from filtered (resp. graded) $\shr$-algebras to filtered 
(resp. graded) $\shr$-stacks, which sends an $\shr$-algebra $\sha$ to the 
stack $\astk\sha$ defined as follows: it is the stack associated with the 
pre-stack $X \supset U\mapsto \astk{\sha(U)}$.
If $f,g\colon \sha\to \shb$ are filtered (resp. graded) $\shr$-algebra 
morphisms, filtered (resp. graded) natural transformations 
$\astk f \Rightarrow \astk g$ are locally described as above.

Let $\sha$ be a filtered $\shr$-algebra. The stack 
$\stkMod[F](\sha)$ of filtered left $\sha$-modules is filtered and equivalent 
to the stack $\stkFun[F](\astk \sha,\stkMod[F](\shr))$ of filtered functors.
Moreover the Yoneda embedding gives a fully faithful functor 
\begin{equation*}
\label{eq:Yoneda}
\astk \sha \to \stkFun[F]((\astk \sha)^\op,\stkMod[F](\shr)) \approx
\stkMod[F](\sha^\op)
\end{equation*}
into the stack of filtered right $\sha$-modules. This identifies $\astk\sha$ 
with the full substack of filtered right $\sha$-modules which are locally 
isomorphic to $\sha$. As above, everything remains true replacing filtered 
algebras and stacks by graded ones.

Let $\stks$ be a filtered $\shr$-stack. We denote by $\stkGr (\stks)$ the 
graded stack associated to the pre-stack $X \supset U\mapsto\catGr {\stks(U)}$.
As above, this defines a functor from filtered $\shr$-stacks to graded ones 
(which is, in fact, a 2-functor). As before, we will make use of the notation
$\stkGr_0$.
\begin{proposition}
Let $\sha$ be a filtered $\shr$-algebra and $\gr (\sha)$ its associated 
graded algebra. Then there is an equivalence of graded 
stacks $\stkGr (\sha^+)\approx\gr (\sha)^+$. 
\end{proposition}

\begin{proof}
Let $\shl$ be a filtered right $\sha$-module locally isomorphic to $\sha$, 
that is, an object of $\sha^+$. Its associated graded module 
$\gr (\shl)$ is a graded right $\gr (\sha)$-module locally isomorphic 
to $\gr (\sha)$, that is, an object of $\gr( \sha)^+$. 
Hence the assignement $\shl \mapsto \gr (\shl)$ induces a functor 
$\stkGr( \sha^+)\to\gr( \sha)^+$ of graded stacks.
Since at each $x\in X$ this reduces to the equality $\catGr {\astk {\sha_x}} 
= \astk{\Gr (\sha_x)}$, it follows that it is a global equivalence. 
\end{proof}

Recall from \cite{Kontsevich2001,D'Agnolo-Polesello2005} that an 
$\shr$-algebroid stack is an $\shr$-stack $\stka$ which is locally non-empty 
and locally connected by isomorphisms. Equivalently, for any  
$x\in X$ there exist an open neighborhood  $U$ of $x$ and an 
$\shr$-algebra $\sha$ on $U$ such that $\stka|_U\approx \astk \sha$. 
Note that, if there exists a global object $L\in\stka(X)$, then 
$\stka\approx\shendo[\stka](L)^+$.
\begin{corollary}
Let $\stka$ be a filtered $\shr$-stack.  If $\stka$ is an $\shr$-algebroid 
stack, then it associated graded stack $\stkGr(\stka)$ is again an 
$\shr$-algebroid stack.
\end{corollary}

Recall that, given an $\shr$-algebroid stack $\stka$, the stack of 
$\stka$-modules is by definition the $\shr$-stack of $\shr$-functors
$\stkFun[\shr](\stka,\stkMod(\shr))$. It is an example of stack of 
twisted sheaves (see for example~\cite{D'Agnolo-Polesello2003,
Kashiwara-Schapira2006}).

\section{Quantization of complex contact manifolds}\label{section:uniq}

Let $(Y,\shl,\alpha)$ be a complex  contact manifold of dimension $2n+1$.
This means that $\shl$ is a line bundle on $Y$ and $\alpha$ is a global 
section of $\Omega_Y^1\tens[\O_Y]\shl$ (an $\shl$-valued 1-form) such that 
$\alpha\wedge (d\alpha)^n$ is a nowhere vanishing global section of 
$\Omega_Y\tens[\O_Y]\shl^{\otimes n+1}$. Recall that a local model for $Y$ is 
an open subset $V$ of $P^*M$, for a complex manifold $M$, equipped with the 
canonical contact structure $(\O_{P^*M}(1)|_V,\lambda)$, where $\lambda$ is 
induced by the Liouville form on $T^*M$. Hence we may define an ($\she,*$)-algebra 
on $Y$ as a filtered sheaf of $\C$-algebras $\sha$ endowed with an isomorphism
of graded algebras $\nu\colon\gr(\sha)\to \DSum_{m\in\Z}\shl^{\otimes m }$ and 
with an anti-involution $\iota$, such that  the triplet $(\sha,\nu,\iota)$ is 
locally isomorphic to $(\opb i\E_M,\sigma,*)$ for any contact local chart 
$i\colon Y\supset V \to P^*M$. 

Although we cannot expect a globally defined ($\she,*$)-algebra on $Y$,
Kashiwara~\cite{Kashiwara1996} proved the existence of a canonical stack 
$\stkMod(\she;Y)$ of microdifferential modules on $Y$. 
Following~\cite{D'Agnolo-Polesello2005}, this result may be restated as:

\begin{theorem}\label{th:existence}
On any complex contact manifold $Y$ there exists a canonical $\C$-stack 
$\stke_Y$ which is locally equivalent to $\astk{(\opb i \E_M)}$ 
for any contact local chart $i\colon Y\supset U \to P^*M$.
\end{theorem}
\noindent
Note that $\stke_Y$ is a $\C$-algebroid stack and $\stkMod(\she;Y)$ is 
equivalent to the stack of $\stke_Y$-modules.

The following proposition allows us to say that the algebroid stack 
$\stke_Y$ provides a quantization of $Y$.
\begin{proposition}
The $\C$-algebroid stack $\stke_Y$ has the following properties:
\begin{itemize}
\item[$\astk {(i)}$] it is filtered;
\item[$\astk {(ii)}$] there is an equivalence of graded stacks 
$$
{\pmb\sigma}\colon\stkGr(\stke_Y)\approxto (\DSum_{m\in\Z}\shl^{\otimes m })^+;
$$
\item[$\astk {(iii)}$] it is endowed with an anti-involution $\pmb *$, that 
is, with an equivalence of stacks
$${\pmb *}\colon \stke_Y^\op\approxto \stke_Y$$
and an invertible transformation $\epsilon\colon {\pmb *}^2 \Rightarrow
\id_{\stke_Y}$ such that the transformations $\epsilon\circ \id_{\pmb *}\colon 
{\pmb *}^3 \Rightarrow {\pmb *}$ and $\id_{\pmb *}\circ\epsilon\colon {\pmb *} 
\Rightarrow {\pmb *}^3$ are inverse one to each other.
\end{itemize}
\end{proposition}
Moreover, $\astk {(i)}$, $\astk {(ii)}$ and $\astk {(iii)}$ are compatible: 
${\pmb *}$ is a filtered $\C$-functor, $\epsilon$ is a filtered natural 
transformation and there is an invertible transformation
$\delta_0\colon \pmb \sigma_0 \circ \stkGr_0 (\pmb *) \Rightarrow 
D\circ \pmb \sigma_0$, where 
$D$ denotes the functor $(\astk {\O_Y})^\op\approxto 
\astk {\O_Y}$, which sends a locally free $\O_Y$-module of rank one to its 
dual, making the following diagram commutative\footnote{The compatibility between 
$\epsilon$ and $\delta_0$ is quite redundant for the algebroid stack $\stke_Y$,
since the data of $(\epsilon,\delta_0)$ are induced by those of $\astk{(\E_M)}$. 
Anyway, it becomes necessary for different choices of $(\epsilon,\delta_0)$.}
\begin{equation}
\label{eq:-1,II}
\xymatrix@C5em{
\pmb \sigma _0 \circ \stkGr_0(\pmb *^2)
\ar@{=>}[r]^-{\delta_0\circ\id_{\stkGr_0(\pmb *)}} 
\ar@{=>}[d]_-{\id_{\pmb \sigma_0}\circ \stkGr_0(\epsilon)} & 
D\circ \pmb \sigma_0\circ \stkGr_0(\pmb *)\\
\pmb \sigma_0 \ar@{=>}[r]^-{\simeq} & D^2 \circ \pmb \sigma_0 
\ar@{=>}[u]_-{\id_{D}\circ \delta_0} .}
\end{equation}

\noindent

\begin{theorem}\label{th:unicity}
The $\C$-algebroid stack $\stke_Y$ is the unique, up to equivalence
-- the equivalence being unique up to a unique isomorphism -- $\C$-stack on $Y$
satisfying the properties $\astk {(i)}$, $\astk {(ii)}$ and $\astk {(iii)}$,
and such that for any contact local chart $i\colon Y\supset V \to P^*M$, 
the triplet $(\stke_Y, {\pmb\sigma},{\pmb *})$ is locally equivalent to 
$((\opb i \E_M)^+,\sigma^+,*^+)$.
\end{theorem}

\begin{proof}
Let $\stkAut[\she,\pmb *](\stke_Y)^{\times}$ be the stack of 
2-groups\footnote{See the Appendix A for the definition, notations and basic 
properties of a stack of 2-groups.} of auto-equivalences of $\stke_Y$ 
preserving the structures  $\astk {(i)}$, $\astk {(ii)}$ and $\astk{(iii)}$.
(Here the upper index $\times$ means that all the non-invertible morphisms 
have been removed.) 
More precisely, its objects are triplets $(\Phi,\beta, \gamma)$, where 
$\Phi\colon \stke_Y\to \stke_Y$ is an equivalence of filtered $\C$-stacks and 
$\beta\colon \pmb \sigma \circ \stkGr (\Phi) \Rightarrow \pmb \sigma$ and
$\gamma\colon  {\pmb *}\circ \Phi \Rightarrow \Phi \circ {\pmb *}$ are
invertible transformations of functors ($\beta$ being graded and $\gamma$ 
filtered) such that the following diagrams commute:
\begin{equation}
\label{eq:gamma}
\xymatrix@C5em{
\pmb * \circ \Phi \circ  {\pmb *}\ar@{=>}[r]^-{\id_{\pmb *}\circ\gamma} 
\ar@{=>}[d]_-{\gamma\circ\id_{\pmb *}} & 
{\pmb *}^2\circ \Phi \ar@{=>}[d]^-{\epsilon\circ\id_{\Phi}} \\
 \Phi\circ {\pmb *}^2 \ar@{=>}[r]_-{\id_{\Phi}\circ \epsilon} &  \Phi}
\end{equation}
\begin{equation}
\label{eq:beta}
\xymatrix@C4em{
\pmb \sigma_0\circ\stkGr_0 (\pmb *)\circ \stkGr_0(\Phi)
\ar@{=>}[rr]^-{\id_{\pmb\sigma_0}\circ \stkGr_0(\gamma)} 
\ar@{=>}[d]_-{\delta_0\circ\id_{\stkGr_0(\Phi)}} && 
\pmb \sigma_0 \circ \stkGr_0(\Phi) \circ  \stkGr_0({\pmb *})
\ar@{=>}[d]^-{\beta \circ \id_{\stkGr_0 (\pmb *)}} \\
D\circ \pmb \sigma_0 \circ \stkGr_0(\Phi) & 
\ar@{=>}[l]^-{\id_{D}\circ \beta} 
D\circ \pmb \sigma_0  & \ar@{=>}[l]^-{\delta_0} 
\pmb \sigma_0 \circ \stkGr_0({\pmb *}).}
\end{equation}
A morphism $\alpha\colon (\Phi_1,\beta_1, \gamma_1)\Rightarrow
(\Phi_2,\beta_2, \gamma_2)$ is a filtered invertible transformation of functors
$\alpha\colon \Phi_1 \Rightarrow \Phi_2$ making the following diagrams 
commutative:
\begin{equation}
\label{eq:alpha}
\xymatrix@C5em{
{\pmb *} \circ \Phi_1 \ar@{=>}[r]^-{\gamma_1}  & 
\Phi_1 \circ \pmb *\ar@{=>}[d]^-{\alpha\circ\id_{\pmb *}} \\
\pmb *\circ\Phi_2 \ar@{=>}[u]^-{\id_{\pmb *}\circ \alpha} 
\ar@{=>}[r]_-{\gamma_2} &  \Phi_2\circ \pmb *}
\qquad
\xymatrix@C3em{
\pmb \sigma \circ \stkGr(\Phi_1)\ar@{=>}[r]^-{\beta_1} 
\ar@{=>}[d]_-{\id_{\pmb\sigma}\circ \stkGr (\alpha)} & \pmb \sigma.\\
\pmb \sigma\circ \stkGr({\Phi_2}) \ar@{=>}[ur]_-{\beta_2} &}
\end{equation}

\noindent
To prove the theorem, it is then enough to show that this stack is trivial, 
that is, consists only in the identity functor $\id_{\stke_Y}$ and in the 
identity transformation $\id_{\stke_Y}\Rightarrow \id_{\stke_Y}$. 

{\em First step.} At any point $y\in Y$, we may find a contact local chart 
$i\colon Y\supset V \to P^*M$ containing $y$ such that $\stke_Y|_V$ is 
equivalent to $(\she^{\sqrt v} _V)^+ = (\opb i\E_M)^+$. We thus reduced to the 
study of the auto-equivalences of $(\she^{\sqrt v} _V)^+$ preserving 
the structures  $\astk {(i)}$, $\astk {(ii)}$ and $\astk {(iii)}$.
Let $(\Phi,\beta, \gamma)$ be such an equivalence. Since the functor 
$(\cdot)^+$ is locally full, up to shrinking the open subset $V$, one may 
suppose that $\Phi$ is isomorphic to $\astk{\varphi}$ for a filtered
$\C$-algebra isomorphism $\varphi\colon \E _V \to \E_V$ and that 
$\pmb\sigma$ (resp. $\pmb *$) is isomorphic to $\astk{\sigma}$ 
(resp. $\astk{*}$), so that $\delta_0$ (resp. $\epsilon$) is the identity 
morphism.
The graded transformation  $\beta$ is thus given by a nowhere vanishing
function $f\in\O_{P^*M}|_V\simeq\gr_0(\E_V)$ and it implies the equalities
$\sigma_m(\varphi(S))=\sigma_m(S)$ for all $S\in \E_V(m)$, so that 
$\gr(\varphi)$ must be the identity. Similarly, the transformation $\gamma$ is given 
by an invertible operator $P \in \E_V$ of order $0$ satisfying
\begin{equation}
\label{eq:P}
\varphi \circ * = \ad(P)\circ *\circ\varphi.
\end{equation}
A direct computation shows that the diagram~\eqref{eq:gamma} corresponds 
to the equality $P^*=P$ and the diagram~\eqref{eq:beta} to $f^2 \sigma_0 (P) =  1$. 

\noindent
By a result of Kashiwara (see~\cite[Lemma 5.3]{Polesello-Schapira} for a 
proof), there exists an invertible operator $Q\in \E_V$ of order $0$
satisfying $Q=Q^*$ and $P=Q^2$. We may choose $Q$ in such a way that 
$f\sigma_0 (Q)=1$. Set $\tilde\varphi= \ad(\opb Q) \circ\varphi$. 
Using~\eqref{eq:P}, one easily gets the relation  
$$
*\circ\tilde\varphi = \tilde\varphi\circ * ,
$$
so that $\tilde\varphi$ is an $\she$-algebra automorphism of $\E_V$.
Moreover, the section $Q$ defines a morphism 
$(\tilde\varphi^+,\id,\id) \Rightarrow (\Phi, \beta,\gamma)$.

\noindent
It follows that the functor of stacks of 2-groups
$$
\left[\she^{\sqrt v, *}_V \to[\ad]  \shaut[\she,*](\she^{\sqrt v}_V) \right]
\longrightarrow \stkAut[\she,\pmb *](\astk{(\she^{\sqrt v}_V)})^\times \qquad
\psi \mapsto (\psi^+,\id,\id)
$$
is locally essentially surjective, where the left-hand side denotes the 
stack of 2-groups associated to the crossed module
$\she^{\sqrt v, *}_V \to[\ad]  \shaut[\she,*](\she^{\sqrt v}_V)$.

{\em Second step.} Let $\psi_1,\psi_2$ be two sections of  
$\shaut[\she,*](\she^{\sqrt v}_V)$ and consider a morphism
$\alpha\colon(\psi_1^+,\id,\id)\Rightarrow (\psi_2^+,\id,\id)$ 
in $\stkAut[\she, \pmb *](\astk{(\she^{\sqrt v}_V)})^\times$. 
Then $\alpha$ is locally given by an invertible operator $R\in\E_V$ of order 
$0$ satisfying
$$
\psi_2 = \ad (R) \circ \psi_1,
$$ 
where the diagrams~\eqref{eq:alpha} correspond to the equalities
$RR^*=1$ and $\sigma_0(R)=1$. Therefore $R$ defines a morphism 
in 
$\left[\she^{\sqrt v, *}_V \to[\ad]  \shaut[\she,*](\she^{\sqrt v}_V) \right]$. 
The above functor is thus locally full. Since it also faithful, we get an equivalence
$$\left[\she^{\sqrt v, *}_V \to[\ad]  \shaut[\she,*](\she^{\sqrt v}_V)\right]\approx 
\stkAut[\she,\pmb *](\astk{(\she^{\sqrt v}_V)})^\times.$$

{\em Third step.} By Lemma~\ref{lemma:key}, the morphism 
$\ad\colon \she^{\sqrt v, *}_V \to \shaut[\she,*](\she^{\sqrt v}_V)$ is an 
isomorphism. Thanks to Proposition~\ref{prop:trivial2-gr}, the associated 
stack of 2-groups is then equivalent to the trivial one. It follows that  
$\stkAut[\she, \pmb *](\stke_Y)^{\times}$ is locally, and hence globally, trivial.
\end{proof}

\begin{corollary}
There exists an ($\she,*$)-algebra on $Y$ if and only if 
$\stke_Y$ has a global object.
\end{corollary}

\begin{proof}
Let $L$ be a global objet of $\stke_Y$, that is, $L\in \stke_Y(Y)$, and set
$\sha=\shendo[\stke_Y](L)$, the sheaf of endomorphisms of $L$. Then
$\sha$ is a filtered $\C$-algebra locally isomorphic to 
$\opb i \she_M$, for any contact local chart  $i\colon Y\supset U \to P^*M$.
By $\astk {(ii)}$, its graded algebra is isomorphic to 
$\DSum_{m\in\Z}\shl^{\otimes m }$. By $\astk {(iii)}$,
there is an isomorphism $*\colon\sha^{\op}\isoto\sha^*$, where 
$\sha^*=\shendo[\stke_Y](L^*)$ is an inner form of $\sha$, that is, 
locally isomorphic to $\sha$ with the glueing automorphisms of the form 
$\ad(P)$ for some $P\in\sha^{\times}$.  Then one may replace $\sha$ by 
a twisted version $\widetilde\sha$, in such a way that $*$ defines an 
anti-involution on it. It follows that $\widetilde\sha$ is an ($\she,*)$-algebra on $Y$.

\noindent
Conversely, let $\sha$ be an ($\she,*$)-algebra on $Y$. Therefore 
$\astk\sha$ is a $\C$-algebroid stack satisfying $\astk {(i)}$,
$\astk {(ii)}$ and $\astk {(iii)}$, and locally isomorphic to 
$\astk{(\opb i \she_M)}$, for any contact local chart 
$i\colon Y\supset U \to P^*M$.
By Theorem~\ref{th:unicity}, it is equivalent to $\stke_Y$. 
Since $\astk\sha$ has a global object, so does $\stke_Y$ .
\end{proof}

\begin{remark}
\begin{itemize}
     \item[(1)] Let $Y= P^*M$. If $\dim M < 3$, there are non-isomorphic 
     ($\she,*$)-algebras (cf~\cite{BoutetdeMonvel1999}), but they all give 
     rise to equivalent algebroid stacks.

     \item[(2)] Using similar techniques, one may show that Theorems~\ref{th:existence} 
     and ~\ref{th:unicity} hold also in the real case, hence replacing $\she$-algebras 
     with Toeplitz algebras. However, in that case the algebroid so obtained has always a 
     global object (cf~\cite{Deligne1995}). It follows that in the real case there 
     are no significative differences between the quantizations given 
     by algebras and those provided by algebroids. 
\end{itemize} 
\end{remark}

\begin{appendix}

\section{Stacks of 2-groups}

In this section we briefly recall the notion of stack of 2-groups, which is, 
roughly speaking, a stack with group-like properties. Reference is made 
to~\cite{Breen1994}\footnote{We prefer to follow here the terminology of Baez-Lauda 
[{\em Higher-dimensional algebra. V. 2-groups}, Theory Appl. Categ. 12 
(2004), 423--491], which seems to us more friendly than the classical one 
of $gr$-category as in {\it loc. cit.}} . We assume that the reader is familiar with the notions 
of monoidal category, monoidal functor and monoidal transformation. 
(The classical reference is~\cite{MacLane} and a more recent one 
is~\cite{Leinster}.)

Let $X$ be a topological space.

\begin{definition}
\begin{itemize}
     \item[(i)] A 2-group
     is a monoidal category $(\mathsf{G},\tens,{\bf 1})$ whose morphisms are 
     all invertible and such that for each object $P\in \mathsf{G}$ there 
     exist an object $Q$ and natural morphisms $P\tens Q\simeq {\bf 1}$ and 
     $Q\tens P\simeq {\bf 1}$. 
     A functor (resp. a natural transformation of functors) of 2-groups is a 
     monoidal functor (resp. a monoidal transformation of functors) between 
     the underlying monoidal categories.
     
     \item[(ii)] A pre-stack (resp. stack) of 2-groups on $X$ is a pre-stack 
     (resp. stack) $\stkg$ such that for each open subset $U\subset X$, the 
     category $\stkg(U)$ is a 2-group, and the restriction functors and the 
     natural transformations between them respect the 2-group structure.
     \end{itemize}
\end{definition}

Let $\shg$ be a sheaf of groups on $X$. We denote by $\shg[0]$ the discrete 
stack defined by trivially enriching $\shg$ with identity arrows, and by $\shg[1]$ 
the stack associated to the pre-stack whose category on an open subset 
$U\subset X$ has a single object $\bullet$  and $\Endo(\bullet)=\shg(U)$ 
as set of morphisms. 
Clearly $\shg[0]$ is a stack of 2-groups, while  $\shg[1]$ defines a stack of 
2-groups if and only if $\shg$ is commutative.

     
\medskip

Recall that a crossed module on $X$ is a complex\footnote{Here we use the convention 
as in~\cite{Breen1994} for which $\shg^{i}$ is in $i$-th degree.} of sheaves of groups 
$\shg^{-1}\to[d]\shg^0$ endowed with a left action of $\shg^0$ on $\shg^{-1}$
such that, for any local sections $g\in\shg^0$ and $h,h'\in\shg^{-1}$, one has
$$
d({}^gh)=\ad(g)(d(h))\quad \text{ and } \quad {}^{d(h')} h=\ad(h')(h).
$$
A morphism of crossed modules is a morphism of complexes compatible with the 
actions in the natural way. 

To each crossed module $\shg^{-1}\to[d]\shg^0$ there is an associated stack of 
2-groups on $X$, which we denote by $[\shg^{-1}\to[d]\shg^0]$, 
defined as follows: it is the stack associated to the pre-stack of 2-groups 
whose objects on an open subset $U\subset X$ are the sections $g\in \shg^0(U)$,
with 2-group law $g\tens g'=gg'$, and whose morphisms $g\to g'$ are given by 
sections $h\in \shg^{-1}(U)$ satisfying $g' = d(h) g$. The 2-group structure 
for morphisms is given by the rule
$$
(g_1\to[h_1] g'_1)\tens (g_2\to[h_2] g'_2)=g_1g_2\to[h_1{}^{g_1}h_2] g'_1g'_2.
$$
In a similar way, a morphism of crossed modules induces a functor of the 
associated stacks of 2-groups.

\begin{remark}
In fact, it is true that any stack of 2-groups on $X$ comes from a crossed 
module. However, this result is not of practical use. We refer to 
\cite{SGA4} for the proof of this fact in the commutative case and to 
\cite{Brown-Spencer} for the non commutative case on $X=\operatorname{pt}$. 
\end{remark}

By definition, an object $\shp$ of $[\shg^{-1}\to[d]\shg^0]$ on 
an open subset $U\subset X$ is described by an open covering $U = 
\Union\nolimits_i U_i$ and sections $\{g_i\}\in \shg^0(U_i)$, subject 
to the relation $g_i = d(h_{ij}) g_j$ on double intersections $U_{ij}$, for 
given sections $\{h_{ij}\}\in \shg^{-1}(U_{ij})$ satisfying $h_{ij}h_{jk}=h_{ik}$ 
on triple intersections $U_{ijk}$.

If $\shg$ is a sheaf of groups on $X$, then the stack of 2-groups defined by the 
crossed module $1 \to \shg$ is naturally identified with $\shg[0]$. 
If moreover $\shg$ is commutative, the complex $\shg \to 1$ is a crossed module 
and the associated stack of 2-groups is identified with $\shg[1]$. 
     
 
\begin{proposition}\label{prop:trivial2-gr}
Let $\shg^{-1}\to[d]\shg^0$ be a crossed module. 
\begin{itemize}
     \item[(1)] If $d$ is surjective, then $[\shg^{-1}\to[d]\shg^0]\approx \ker d[1]$ as stacks of 2-groups.
     \item[(2)] If $d$ is injective, then $[\shg^{-1}\to[d]\shg^0]\approx \coker d[0]$ as stacks of 2-groups.
\end{itemize}
In particular, if $d$ is an isomorphism, then $[\shg^{-1}\to[d]\shg^0]$ is trivial, {\em i.e.} it is 
equivalent to the stack with one object and one morphism. 
\end{proposition}

\begin{proof}
(1) Consider the fully faithful functor of stacks of 2-groups
$\ker d[1] \to[] [\shg^{-1}\to [d]\shg^0]$,
which sends $\bullet$ to the identity object ${\bf 1} = 1 \in \shg^0$.
Let us check that is is locally essentially surjective. 
Let $g\in \shg^0$ be a local section. We may suppose that there exists a local 
section $h\in \shg^{-1}$ satisfiyng $d(h)=g$. Hence $h$ defines a morphism 
$1\to g$ in $[\shg^{-1}\to[d]\shg^0]$. It follows that any other object of
$[\shg^{-1}\to[d]\shg^0]$ is locally isomorphic to ${\bf 1}$.

(2) Consider the locally essentialy surjective functor of stacks of 2-groups
$[\shg^{-1}\to [d]\shg^0]\to \coker d[0]$, 
which sends an objet to its isomorphism class. It remains to check that 
$[\shg^{-1}\to [d]\shg^0]$ has only trivial arrows. But this is clear, since 
for any object there are no automorphisms other than the identity 
$1\in\shg^{-1}$.
\end{proof}

\end{appendix}

\providecommand{\bysame}{\leavevmode\hbox to3em{\hrulefill}\thinspace}

\end{document}